\documentstyle[12pt]{article} \vsize=28.7truecm \hsize=23truecm \columnsep=0.8truecm
 \def\vbar{\mathchoice{\vrule height2.3ptdepth-.3ptwidth.12pt\kern-
 .10pt}
    {\vrule height6.3ptdepth-.3ptwidth.11pt\kern-.11pt}
    {\vrule height5.1ptdepth-.30ptwidth.8pt\kern-.8pt}
    {\vrule height4.1ptdepth-.24ptwidth.6pt\kern-.7pt}}
\def\reel{\hbox{I\hskip-2pt R}}

\def\<{\langle}
\def\>{\rangle}

\def\n{{\boldmath n}}

\def\<{\langle}
\def\>{\rangle}

\def\reel{\hbox{I\hskip -2pt R}}

\def\n{{\noindent}}
\newtheorem{theorem}{Theorem}[section]
\newtheorem{definition}{Definition}[section]

\newtheorem{lemma}[theorem]{Lemma}
\newtheorem{prop}[theorem]{Proposition}

\newtheorem{remark}{Remark}[section]
\newtheorem{example}{Example}[section]

\newenvironment{Pff}{\hspace*{-\parindent}{\bf Proof}}{\hfill $\Box$
\vspace*{0.2cm}}

\def\opp{\;\mathrm{opp}\;}

\begin{document}
 \vspace{0.5cm}
\title{{\bf
Riesz exponential families on homogeneous cones}}
\author{I. Boutouria$^*$, A. Hassairi\thanks{University of Sfax, Laboratory of Probability and Statistics, B.P. 1171, Sfax, Tunisia}
 \footnote{Corresponding author.
 \textit{E-mail address: abdelhamid.hassairi@fss.rnu.tn}}$\;$
}
 \date{}

 \maketitle

\n {\small \textbf{Abstract.}  In this paper, we introduce, for a
multiplier $\chi$, a notion of generalized power function $x\mapsto
\Delta _{\chi }(x),$ defined on the homogeneous cone ${\mathcal{P}}$
of a Vinberg algebra ${\mathcal{A}}$. We then extend to
${\mathcal{A}}$ the famous Gindikin result, that is we determine the
set of multipliers $\chi$ such that the map $\theta \mapsto \Delta
_{\chi }(\theta ^{-1})$, defined on ${\mathcal{P}}^{\ast}$, is the
Laplace transform of a positive measure $R_{\chi }$. We also
determine the set of $\chi $ such that
$R_{\chi }$ generates an exponential family, and we calculate the variance function of this family.\\
}
 \\\\ Key
words: Homogeneous cone; multiplier; generalized power; Riesz probability distribution; exponential family; variance function.\\

\section{Introduction}
It is well known (see Casalis and Letac (1996)) that the Wishart
distributions on the cone of $(r,r)$ positive symmetric matrices
or on the symmetric cone $\Omega $ of any Euclidean Jordan algebra
$E$ of rank $r$ are the elements of the natural exponential
families generated by the measures $\mu_{p}$ such that the Laplace
transform is defined on $\Omega$ by
$$L_{\mu_{p}}(\theta)=(\det(\theta^{-1}))^{p},$$ for $p$ in
$\displaystyle\{\frac{1}{2},1,\frac{3}{2},\cdot\cdot\cdot,\frac{r-1}{2}\}\cup\displaystyle]\frac{r-1}{2},
\ +\infty[.$ The measure $\mu_{p}$ is absolutely continuous when
$p\in \displaystyle]\frac{r-1}{2}, \ +\infty[$ and is singular
concentrated on the boundary of the cone, when
$p\in\displaystyle\{\frac{1}{2},1,\frac{3}{2},\cdot\cdot\cdot,\frac{r-1}{2}\}$.
In 2001, Hassairi and Lajmi have introduced the Riesz distribution
on $\Omega $ as an extension of the Wishart distribution. These
authors have started from the fact that in a Jordan algebra,
besides the real power of the determinant, there is the so called
generalized power $\Delta _{s}(x)$ of an element $x$ of $\Omega ,$
defined for a fixed ordered Jordan frame of $E$ and for
$s=(s_1,\cdot\cdot\cdot,s_r)$ in ${\reel}^r$, and they have used a
remarkable result, due to Gindikin (1964), which determines the
set $\Xi $ of $s$ in ${\reel}^r$ such that $\Delta _s(\theta
^{-1})$ is the Laplace transform of some positive measure $R_s$ on
$E$. The generalized power $\Delta _{s}(x)$ is
a power function of the principal minors of $x$ which reduces to $(\det (x))^p$ in the particular case where $%
s_1=s_2=\cdot\cdot\cdot=s_r=p$, and in this case, the measure
$R_s$ in nothing but $\mu_{p}$. We mention here that Ishi (2000)
has given a more detailed description of the Gindikin set $\Xi$
based on the orbit structure of $\overline{\Omega }$ under the
action of some Lie group. He has also given explicitly the measure
$R_s$ for each $s$ in $\Xi$. In all these works, the definition of
the Riesz measure $R_s$ and in particular of the Riesz probability
distribution is based on the choice of a totally ordered Jordan
frame which allows the definition of the principal minors and of
the generalized power of an element of the algebra. The fact that
the order is total is a fundamental condition not only for the
definition of the distribution but also in the proof of many
results. To define models in which some specified conditional
independencies, usually given by a graph, are taken into account,
there has been an interest in probability distributions on the
homogeneous cone of a Vinberg algebra. For instance, Andersson and
Wojnar (2004) have defined a class of absolutely continuous
$``$Wishart" distributions on an homogeneous cone. These
distributions have been characterized by Boutouria (2005, 2007) in
the Bobecka and Weso\l owski (2002) way. They have also been
charcterized by Boutouria and Hassairi (2008) in the way given in
Olkin and Rubin (1962) for the ordinary Wishart. The aim of the
present work is to use an approach similar to the one used in the
definition of a Riesz exponential family on a symmetric cone to
introduce a Riesz exponential family on an homogeneous cone. The
distributions in these families are defined for any graph, that is
for any order relation not necessary total. Some of these
distributions are absolutely continuous with respect to the
Lebesgue measure and some are singular concentrated on the
boundary of the cone. In this connection, the Riesz distribution
on the symmetric cone of a Jordan algebra
may be seen as the particular one corresponding to the particular directed graph with vertex set $%
\left\{ 1,\cdot\cdot\cdot,r\right\}$ and edges defined by the usual
order on integers. We first define for  an element of an homogeneous
cone two kinds of principal minors, minors which are said strict and
minors which are said large. We then define for a multiplier $\chi
$, a notion of generalized power function $x\mapsto \Delta _{\chi
}(x).\;
$One of our main results is the determination of the set of multipliers $%
\chi $ such that the map $\theta \mapsto \Delta _{\chi }(\theta
^{-1})$ is the Laplace transform of a positive measure $R_{\chi }$.
It is a generalization of  Gindikin result with a more elaborate
proof adapted to the properties of the Vinberg algebra and the
graph. Concerning the generated exponential families, we give a
necessary and sufficient condition on $\chi $ in order that $R_{\chi
} $ generates an exponential family and, under this condition, we
determine the variance function of the family.
\section{Vinberg algebras and homogeneous cones}
In this section, we introduce some notations and review some basis
concepts concerning Vinberg algebras and their homogeneous cones. We
also introduce a useful decomposition of an element of the cone.

Throughout the paper, $I$ denote a partially ordered finite set
equipped with a relation denoted $\preceq$. We will write $i\prec j$
if $i\preceq j$ and $i \not = j$. For all pairs $(i,j)\in I \times
I$ with $j\prec i$, let $E_{ij}$ be a finite-dimensional vector
space over $\reel$ with $n_{ij}=\mbox{dim}(E_{ij})>0$. Set
\begin{eqnarray*}
{\mathcal{A}}_{ij}= \left \{ \begin{array}{ccccc}
  {\reel} & \textrm{for} & i=j&& \\
  E_{ij} & \textrm{for} & j\succ i & \textrm{or} & j\prec i\\
  \{0\} &  &\textrm{otherwise}.&&
\end{array}\right.
\end{eqnarray*}

\n and ${\mathcal{A}}=\displaystyle\prod\limits _{i,j \in I\times
I}{\mathcal{A}}_{ij}\cdot$ An element $A\equiv (a_{ij}, \ i,j \in
I)$ of ${\mathcal{A}}$ may be seen as a matrix and so we define
the trace $\textrm{tr} A=\displaystyle \sum _{i \in I}a_{ii}$. We
also define\begin{eqnarray}\label{les n_i}n_{i.}=\displaystyle\sum
\limits _{\mu < i} n_{i\mu}, n_{.i}=\displaystyle\sum \limits _{i<
\mu }n_{\mu i}, \ n_{i}=1+\frac{1}{2}(n_{i.}+ n_{.i}), i\in I \
\textrm{and}\ n_{.}=\displaystyle\sum \limits _{i \in
I}n_{i}\cdot\end{eqnarray}\n Let $f_{ij}\;:\; E_{ij} \rightarrow
E_{ij}$, $i\succ j$, be
 involutional linear mappings, i.e., $f_{ij}^{-1}=f_{ij}$. They
induce an involutional mapping ( $ A \mapsto A^{\ast}$ ) of
${\mathcal{A}}$ given as follows: $A^{\ast}=(a^*_{ij}|(i,j)\in I
\times I)$, where

$$a^*_{ij}=\left \{\begin{array}{ccccc}a_{ii} \hskip1.5cm& \textrm{for} & i=j&&  \\
 f_{ij}(a_{ij})=a_{ij}^{\ast} & \textrm{for} &j\prec i&\textrm{or}&i\prec j \\
  \{0\}\hskip1.5cm &  & \textrm{otherwise}.&& \end{array}\right.$$

\n Let ${\mathcal{T}}_{u}=\{A\equiv(a_{ij})\in{\mathcal{A}},\
\forall i,j \in I :i \not\preceq j \Rightarrow a_{ij}=0\},$
${\mathcal{T}}_{l}=\{A\equiv(a_{ij})\in{\mathcal{A}},\ \forall i,j
\in I :j\not\preceq i  \Rightarrow a_{ij}=0\}$ and $\mathcal{H}=\{A
\in{\mathcal{A}}, \ A^{\ast}=A\}$ denote respectively the set of
upper triangular matrices, the lower triangular matrices and the
Hermitian matrices.
 \n The sets of upper and lower triangular matrices in
${\mathcal{P}}$ with positive diagonal elements are respectively
denoted by ${\mathcal{T}}_{u}^{+}$ and ${\mathcal{T}}_{l}^{+}$. The
sets of diagonal matrices and of diagonal matrices with positive
entries are denoted by ${\mathcal{D}}$ and ${\mathcal{D}}^{+}$,
respectively.

The space ${\mathcal{A}}$ is equipped with a bilinear map called
multiplication and denoted by $(A,B)\mapsto AB$, using bilinear
mappings ${\mathcal{A}}_{ij}\times {\mathcal{A}}_{jk}\rightarrow
{\mathcal{A}}_{ik}$, denoted by $(a_{ij},b_{jk})\mapsto
a_{ij}b_{jk}, $ such that $AB=C\equiv (c_{ij}|(i,j)\in I\times I)$
with $c_{ij}=\displaystyle\sum \limits _{\mu \in I}a_{i\mu}b_{\mu
j}$.

\n The multiplication is required to satisfy the following
properties:
\begin{eqnarray*}
\label{vinberg} &&i)\ \  \forall A\in {\mathcal{A}}; \ A \not= 0
\Rightarrow
\textrm{tr}(AA^{\ast})>0\nonumber\\
&&ii) \ \  \forall A, B \in {\mathcal{A}}; \
(AB)^{\ast}=B^{\ast}A^{\ast}\nonumber\\
&&iii)\ \   \forall A, B \in {\mathcal{A}}; \ \textrm{tr} (AB)=\textrm{tr}(BA)\\
&&iv) \ \  \forall A, B, C \in {\mathcal{A}}; \  \textrm{tr} (A(BC))=\textrm{tr} ((AB)C)\nonumber\\
&&v)\ \   \forall U,S,T \in {\mathcal{T}}_{l};  \ (ST)U=S(TU)\nonumber\\
&&vi) \ \  \forall U,S,T \in {\mathcal{T}}_{l}; \
T(UU^{\ast})=(TU)U^{\ast}.\nonumber
\end{eqnarray*}
An algebra ${\mathcal{A}}$ with the above structure and properties
is called a Vinberg algebra (For more details, we can refer to
Andersson and Wojnar (2004). Define the inner products
$(.,.)_{ij}$ on $E_{ij}$, $i\succ j$ by $\parallel
a_{ij}\parallel^{2}_{ij}=a_{ij}f_{ij}(a_{ij})$, $a_{ij} \in
E_{ij}$. Thus instead of specifying the bilinear form $(a_{ij},
b_{ji})\mapsto a_{ij}b_{ji}$ on $E_{ij}$ one can specify an inner
product $(.,.)_{ij}$ on $E_{ij}$, $i\succ j$. It can be
established that the following two conditions also must hold:
\begin{enumerate}
\item $\forall \ a_{ij} \in E_{ij}, \ b_{jk} \in E_{jk}: \
\parallel a_{ij}b_{jk}\parallel_{ik}^{2}=\parallel
a_{ij}\parallel^{2}_{ij}\parallel b_{jk}\parallel_{jk}^{2}, \ k\prec
j\prec i,$

\n and \item If $a_{ik} \in E_{ik}, \ b_{jk} \in E_{jk}$, with
$k\prec j\prec i$ and $(a_{ik}, c_{ij}b_{jk})_{ik}=0$ for all
$c_{ij} \in E_{ij},$ then $(d_{li}a_{ik}, c_{lj}b_{jk})_{lk}=0$ for
all $l \in I$ with $i\prec l, \mbox{and all} \ c_{lj} \in E_{lj}, \
\mbox{and} \ d_{li} \in E_{li}$.
\end{enumerate}

\n We consider the element $(a_{ij}|(i,j)\in I\times I)$ of
${\mathcal{D}}$ such that $a_{ii}=1, \ \forall \ i \in I$ as the
unit element of ${\mathcal{A}}$ and we denote it by $e$. We also
define $E_k=(d_{ij})_{i, j \in I}\in {\mathcal{D}}$ with $d_{kk}=1$
and $d_{jj}=0 \ \forall j\not=k$. It is clear that
$\displaystyle\sum_{k \in I}E_k=e$.

Vinberg (1965) proved that the subset
${\mathcal{P}}=\{TT^{\ast}\in {\mathcal{A}}, \ T \in
{\mathcal{T}}_{l}^{+}\}$ $\subset \mathcal{H} \subset \mathcal{A}$
forms a homogeneous cone, that is the action of its automorphism
group is transitive. Let $G$ be the connected component of the
identity in Aut$({\cal{P}})$; the group of linear transformations
leaving ${\cal P}$ invariant. We recall that $\chi: G\mapsto
\reel_+$ is said to be a multiplier on the group $G$ if it is
continuous, $\chi(e)=1$ and $\chi(g_1g_2)=\chi(g_1)\chi(g_2)$ for
all $g_1,g_2 \in G$. Consider the map $\pi\;:\; T\in {\mathcal
T}_{l}^+ \mapsto \pi(T)\in \pi({\mathcal T}_{l}^+)\subset G$ such
that for $X=VV^{\ast}\in {\mathcal{P}}$, $V\in {\mathcal T}_{l}^+$
\begin{eqnarray}\label{pi}
\pi(T)(X)=(TV)(V^*T^*).
\end{eqnarray}
Andersson and Wojnar (2004) have shown that the restriction of a
multiplier $\chi$ to the (lower) triangular group ${\mathcal
T}_{l}^+$, i.e., $\chi \circ \pi : T_{l}^{+} \rightarrow
{\reel}_{+}$ is in one to one correspondence with the set of
$(\lambda_{i}, i \in I)\in {\reel}^{I}$. We will then describe a
multiplier $\chi$ by its corresponding point in $\reel^{I}$ and we
denote ${\cal{X}}=\{\chi : T_{l}^{+} \rightarrow {\reel}_{+}\}$.

If $\preceq^{\opp}$ is the opposite ordering on the index set $I$,
i.e., $i\preceq^{\opp}j\Leftrightarrow j\preceq i$. The Vinberg
algebra ${\cal{A}}^{\opp}=\displaystyle\prod\limits _{i,j \in
I\times I}{\mathcal{A}}^{\opp}_{ij}$, where
\begin{eqnarray*}
{\mathcal{A}}^{\opp}_{ij}= \left \{ \begin{array}{ccccc}
  {\reel} & \textrm{for} & i=j&& \\
  E_{ij} & \textrm{for} & j\succ^{\opp} i & \textrm{or} & j\prec^{\opp} i\\
  \{0\} &  &\textrm{otherwise,}&&
\end{array}\right.
\end{eqnarray*}
differs from the Vinberg algebra ${\mathcal{A}}$ only in the
ordering of the index set $I$. Vinberg (1965) proved that
${\mathcal{P}}_{\preceq^{\opp}}=\{T^{\ast}T\in {\mathcal{A}}, \ T
\in {\mathcal{T}}_{l}^{+}\}$ is the dual cone of ${\cal{P}}$. The
inner product $(A,B) \rightarrow \textrm{tr}(AB)$ on $H$
identifies $H$ with its dual $H^{\ast}$, i.e.,
\begin{eqnarray*}
H &\leftrightarrow& H^{\ast}\\
A &\mapsto&(B \mapsto \textrm{tr}(AB)),
\end{eqnarray*}
and this isomorphism identifies ${\mathcal{P}}_{\preceq^{\opp}}$
with the dual cone ${\mathcal{P}}^{\ast}$ of ${\mathcal{P}}$.

Now, for $i \in I$, we denote $I_{\preceq i }=\{j \in I; \ j \preceq
i\}$ and $I_{i \preceq }=\{j \in I; \ i \preceq j\}$ and we say that
$j$ separates $i_{1}$ and $i_{2}$ if $j \in I_{i_{1}\preceq } \cap
I_{i_{2}\preceq }$ and $j\not\in \{i_{1}, \ i_{2}\}$. In this case,
$j$ is called a separator. We denote $S_{i}=\{j \in I_{i \preceq };
\ j \ \textrm{is a separator}\}$ and
${\mathcal{S}}=\displaystyle\bigcup_{i \in I}S_{i}$.

\n If $T=(t_{ij})_{i,j \in I}$ is in ${{\cal T}}_{l}$, we define the
element $T_{i\preceq }$ of ${{\cal T}}_{l}$ by
\begin{eqnarray}\label{T_i}T_{i\preceq } =(t'_{ij})_{i,j \in I}, \
\textrm{with} \ t'_{jk}=t_{jk}\ \textrm{if} \ i\preceq j, \ k  \
\textrm{and} \ t'_{jk}=0 \ \textrm{otherwise},\end{eqnarray} and the
element $T_{i\prec }$ of ${{\cal T}}_{l}$ by
\begin{eqnarray}\label{T_î}T_{i\prec } =(t'_{ij})_{i,j \in I}, \
\textrm{with} \ t'_{jk}=t_{jk}\ \textrm{if} \ i\prec j, \ k  \
\textrm{and} \ t'_{jk}=0 \ \textrm{otherwise}.\end{eqnarray} If
$X=TT^{\ast}$, we denote
\begin{eqnarray}\label{Z_iprece}
X_{i\preceq }=T_{i\preceq }T_{i\preceq }^{\ast},\hskip2cm X_{i\prec
}=T_{i\prec }T_{i\prec }^{\ast}\cdot\end{eqnarray} We also denote by
${\cal{P}}_{i \preceq}$ (resp ${\cal{P}}_{i \prec}$) the set of $
X_{i\preceq }=T_{i\preceq }T_{i\preceq }^{\ast}$ (resp $X_{i\prec
}=T_{i\prec }T_{i\prec }^{\ast})$ corresponding to $T \in {{\cal
T}}_{l}^+$. It is easy to see that ${\cal{P}}_{i \preceq}$ and
${\cal{P}}_{i \prec}$ are respectively the homogeneous cones of the
Vinberg subalgebras of ${\cal A}$ defined by ${\cal{A}}_{i
\preceq}=\displaystyle \prod_{k,l \in I_{i\preceq }}{\cal A}_{kl}$
and ${\cal{A}}_{i \prec}=\displaystyle \prod_{k,l \in I_{i\prec
}}{\cal A}_{kl}$. We denote by $e_{i}$ and $\check{e}_{\check{i}}$
respectively, the unit element of ${\cal{A}}_{i\preceq}$ and
${\cal{A}}_{i\prec}$. We also define the rank of ${\cal{P}}_{i
\preceq}$ (resp the rank of ${\cal{P}}_{i \prec}$) the cardinal of
the set $\{j \in I, i\preceq j\}$ (resp the cardinal of the set $\{j
\in I, i\prec j\}$). Finally, if we denote $\wp=\{i \in I; \ I_{i
\prec}\not=\emptyset \ \textrm{and}\ I_{ \prec i}=\emptyset\}$ and
if we set, for $X \in {\cal{P}}$,
\begin{eqnarray}\label{X_i}
{X}_{i}= \left \{ \begin{array}{ccccc} X_{i\preceq }-
\displaystyle\sum_{s \in S_{i}} X_{s\preceq }& \textrm{if} \ & i \in
\wp
\hskip1.8cm&&\\X_{i\preceq }\hskip1.9cm& \textrm{if}\ &i \in S\hskip1.8cm&& \\
{0} \hskip2.1cm&& \textrm{otherwise},\hskip0.9cm && \\
\end{array}\right.
\end{eqnarray}
\n then, we have the following decomposition of $X$
\begin{eqnarray}\label{decomposition}
X=\displaystyle\sum_{i \in I}X_{i}.
\end{eqnarray}

\section{Riesz measures on an homogeneous
cone} The definition of a Riesz measure on the symmetric cone of a
Jordan algebra relies on the notion of generalized power of an
element of the cone which is a power function of the so-called
principal minors. In order to define a Riesz distribution on an
homogeneous cone, we need to extend all these things to a Vinberg
algebra.
\subsection{Generalized power}
We first introduce a notion of determinant. For $X=TT^\ast$, with
$T=(t_{ij}) \in {\mathcal{T}}_{l}$, we define the determinants
$$\mbox{det} X= \displaystyle \prod _{i \in I} t_{ii}^{2}, \hskip
1cm \mbox{det}^{\preceq}X_{\preceq i}=\displaystyle \prod _{j\in
I_{\preceq i}} t_{jj}^{2}\  \hskip 0.5cm\textrm{and} \hskip 0.5cm
\mbox{det}^{\prec}X_{\prec i}=\displaystyle \prod _{j \in I_{\prec
i}} t_{jj}^{2}.$$

\n For $X=TT^\ast$, with $T \in {\mathcal{T}}_{l}$, we define the
strict principal minor of order $k$ of $X$  as
\begin{eqnarray}\label{deta prec k}
\Delta_{\prec k}(X)= \left \{ \begin{array}{ccccc}
   \det^{\prec}(X_{\prec k})& \textrm{if} & I_{\prec k}\not=\emptyset&& \\
  1 & \textrm{if} & I_{\prec k}=\emptyset\\
\end{array}\right.
\end{eqnarray}
and the large principal minor of  order $k$ of $X$ as
\begin{eqnarray}\label{deta preceq k}
\Delta_{\preceq k}(X)=\mbox{det}^{\preceq}(X_{\preceq k}).
\end{eqnarray}

\begin{definition} Let $\chi=\{\lambda_i, i \in I\}$ be a multiplier and $X \in {\cal{P}}$, then the map defined by
\begin{eqnarray}\label{genrateur power}
X\mapsto\Delta_{\chi}(X)=\displaystyle\prod_{k \in
I}(\frac{\Delta_{\preceq k}(X)}{\Delta_{\prec k}(X)})^{\lambda_{k}}.
\end{eqnarray}
is called the generalized power function corresponding of $\chi$.
\end{definition}
We also denote
\begin{eqnarray}\label{genra power p_i}
\Delta_{\chi}^{(i)}(X)=\displaystyle\prod_{k \in I_{i\preceq
}}(\frac{\Delta_{\preceq k}(X)}{\Delta_{\prec k}(X)})^{\lambda_{k}}.
\end{eqnarray}
Note that, if $\lambda_{i}=\lambda, \ i \in I$,  then
$\Delta_{\chi}(X)=(\det X)^{\lambda}$. It is easy to verify that
$\Delta_{\chi+\chi'}(X)=\Delta_{\chi}(X)\Delta_{\chi'}(X)$, where
$\chi+\chi'=\{\lambda_{i}+\lambda'_{i}, \ i \in I\}\cdot$
\begin{example}
Let us consider $I=\{1,2,3,4\}$ and the  poset defined by
$$1\prec 3, \ 1 \prec 4, \ 2\prec 3.$$
For $X =TT^{\ast}\in {\cal{P}}$, with $T=(t_{ij}) \in
{\mathcal{T}}_{l}$, we have $\Delta_{\preceq
1}(X)=\mbox{det}^{\preceq}(X_{\preceq 1})=t_{11}^{2}$,
$\Delta_{\preceq 2}(X)=\mbox{det}^{\preceq}(X_{\preceq
2})=t_{22}^{2}$, $\Delta_{\preceq
4}(X)=\mbox{det}^{\preceq}(X_{\preceq 4})=t_{11}^{2}t_{44}^{2}$ and
$\Delta_{\preceq 3}(X)=\mbox{det}^{\preceq}(X_{\preceq
3})=t_{11}^{2}t_{22}^{2}t_{33}^{2}$. Hence, for $\chi=\{\lambda_i=1,
\ i \in I\}$,
\begin{eqnarray*}\Delta_{\chi}(X)=\frac{t_{11}^{2}}{1}\frac{t_{22}^{2}}{1}\frac{t_{11}^{2}t_{44}^{2}}{t_{11}^{2}}\frac{t_{11}^{2}t_{22}^{2}t_{33}^{2}}{t_{11}^{2}t_{22}^{2}}=t_{11}^{2}t_{22}^{2}t_{33}^{2}t_{44}^{2}.
\end{eqnarray*}
\end{example}
\subsection{Orbit decomposition of the closure
$\overline{{\cal{P}}}$ of ${\cal{P}}$}

For $i \in \wp\cup \cal{S}$, we denote by $\varepsilon^{i}$ the
set of maps $\psi$ defined from $I$ into $\{0, 1 \}$ as follows:
If $i \in \wp$, $\psi$ is such that $\psi(j)=0$, when $j \not \in
I_{i\preceq}$ or $j \in \cal{S}$, and if $i \in \cal{S}$, $\psi$
is such that $\psi(j)=0,$ when $j \not \in I_{i\preceq}$.
Similarly, we denote by $\varepsilon^{\check{i}}$ the set of maps
$\psi$ defined from $I$ into $\{0, 1 \}$. If $i \in \wp$, $\psi$
is such that $\psi(j)=0$, when $j \not \in I_{i\prec}$ or $j \in
\cal{S}$, and if $i \in \cal{S}$, $\psi$ is such that $\psi(j)=0,$
when $j \not \in I_{i\prec}$. With these notations, we define for
$i \in \wp\cup \cal{S}$ and $\psi \in \varepsilon^{i}$, $e_{\psi}
=\textrm{diag}(\psi)=\textrm{diag}(\psi(j), j \in I)\in {\cal{D}}$
and we denote by
${\mathcal{T}}_{l}^{+}.e_{\psi}=\{Te_{\psi}T^{\ast}, \ T \in
{\mathcal{T}}_{l}^{+}\}$. We also consider the two elements of
$\cal{A}$
$${E}^{i}= \left\{\begin{array}{ccccc}
  e_{i} & \textrm{in} \hskip0.5cm {\cal{A}}_{i \preceq}&& \\
  0 &\textrm{elsewhere },&&\\
\end{array}\right. \textrm{and}\hskip0.5cm \  \check{{E}}^{\check{i}}= \left\{\begin{array}{ccccc}
 \check{ e}_{\check{i}} & \textrm{in} \hskip0.5cm {\cal{A}}_{i \prec}&& \\
  0 &\textrm{elsewhere }.&&\\
\end{array}\right.$$
Next, we state and prove a fundamental result. It is a decomposition
of $\overline{{\cal{P}}}$ in orbits.
\begin{theorem} $i)$\begin{equation}\label{p bar}
\overline{{\cal{P}}}=\displaystyle\sum_{i \in \wp\cup
\cal{S}}\overline{{\cal{P}}}_{i \preceq}\cdot
\end{equation}
$ii) $ Let $i \in \wp\cup \cal{S}$, then
\begin{equation}\label{ipns p_i bar}
\overline{{\cal{P}}}_{i \preceq}=\displaystyle\bigcup_{\psi \in
\varepsilon^i }{\mathcal{T}}_{l}^{+}.e_{\psi}.
\end{equation}
\end{theorem}

\begin{Pff} \ $i)$ Let $Z \in \overline{{\cal{P}}}$, then there exist a
sequence $\{Z^{(n)}\}_{n \in N } $ in ${\cal{P}}$ such that
$Z^{(n)}\rightarrow Z$ as  $n\rightarrow \infty$. Since
$\{Z^{(n)}\}_{n \in N } $ is in ${\cal{P}}$, using the decomposition
(\ref{decomposition}) we write $Z^{(n)}=\displaystyle\sum_{i\in
I}Z_i^{(n)}$, where $Z_i^{(n)}\in \overline{{\cal{P}}}_{i\preceq}$.
Hence $Z=\displaystyle\sum_{i \in I}Z_i$, where $Z_i \in
\overline{{\cal{P}}}_{i\preceq}$ and (\ref{p bar}) is proved.

\n $ii)$ We will prove (\ref{ipns p_i bar}) by induction on the rank
of the cone ${\cal{P}}_{i \preceq}$. It is obvious that (\ref{ipns
p_i bar}) holds for $i\in \wp\cup \cal{S}$ such that
$\mbox{rank}{\cal{P}}_{i\preceq}=1$. Suppose that (\ref{ipns p_i
bar}) holds for any $i\in \wp\cup \cal{S}$ such that
$\mbox{rank}{\cal{P}}_{i\preceq}<l$ and let us show that it holds
for $i$ such that $\mbox{rank}{\cal{P}}_{i\preceq}=l$. Consider the
set $M_i=\{j\in I; \ I_{\prec j}=\{i\}\}$. Then using the
decomposition defined by (\ref{X_i}) and (\ref{decomposition}) for
an element of the cone ${\cal{P}}_{i \prec}$, we easily see
\begin{eqnarray}\label{P_M}
\overline{{\cal{P}}}_{i \prec}=\displaystyle\sum_{j \in
M_i}\overline{{\cal{P}}}_{j \preceq}\cdot \end{eqnarray} As
$\mbox{rank}{\cal{P}}_{i\preceq}=l$, we have that
$\mbox{rank}{\cal{P}}_{i\prec}=l-1$ and it follows that
$\mbox{rank}{\cal{P}}_{j\preceq}\leq l-1$, $\forall j \in M_i$.
Using the induction hypothesis, we can write \begin{eqnarray*}
\overline{{\cal{P}}}_{j \preceq}=\displaystyle\bigcup_{\psi \in
\varepsilon^j }{\mathcal{T}}_{l}^{+}.e_{\psi}, \ \hskip1cm j\in
M_i\cdot
\end{eqnarray*}Now, let $\varepsilon^{i\prec}=\displaystyle\sum_{j \in
M_i}\varepsilon^{j}$, then we obtain
\begin{eqnarray}\label{p_iprec}
\overline{{\cal{P}}}_{i \prec}=\displaystyle\sum_{j \in
M_i}\overline{{\cal{P}}}_{j \preceq}=\displaystyle\sum_{j \in
M_{i}}(\displaystyle\bigcup_{\psi \in \varepsilon^j
}{\mathcal{T}}_{l}^{+}.e_{\psi})=\displaystyle\bigcup_{\psi \in
\varepsilon^{i\prec}
}{\mathcal{T}}_{l}^{+}.e_{\psi}\subset\displaystyle\bigcup_{\psi \in
\varepsilon^{i}} {\mathcal{T}}_{l}^{+}.e_{\psi}\cdot
\end{eqnarray}To conclude, we will verify that for $Z \in \overline{{\cal{P}}}_{i
\preceq}$, there exist $\psi\in \varepsilon^i$ and $T \in
{\mathcal{T}}_{l}^{+}$ such that $Z=T.e_{\psi}$. Let $Z \in
\overline{{\cal{P}}}_{i \preceq}$, then there exists a sequence
$\{Z^{(n)}\}_{n\in N } $ in ${\cal{P}}_{i \preceq}$ such that
$Z^{(n)}\rightarrow Z$ as $n\rightarrow \infty$. As $Z^{(n)} \in
{\cal{P}}_{i \preceq}$, there exists $U_{i\preceq
}^{(n)}=(u_{jk}^{(n)})_{j,k \in I}$ in ${\cal{T}}_{l}$ such that
$Z^{(n)}=U_{i\preceq }^{(n)}(U_{i\preceq }^{(n)})^{\ast}$ (see
(\ref{T_i}) and (\ref{Z_iprece})). In particular, we have
\begin{equation}
z^{(n)}_{kk}=(u_{kk}^{(n)})^{2}+\displaystyle\sum_{j\prec
k}\parallel u^{(n)}_{kj}\parallel_{kj}^{2},
\end{equation}
for $k \in I_{i \preceq}$. This implies that the sequences
$(u_{kk}^{(n)})_{n \in {\bf N}}$ and $(u_{kj}^{(n)})_{n \in {\bf
N}}$ are bounded. Therefore there exists a subsequence of positive
integers $(n_m)$  such that $(u_{kk}^{(n_m)})_{m}$ and
$(u_{kj}^{(n_m)})_{m}$ converge. Let
$\widetilde{u}_{kk}=\displaystyle\lim_{ m\rightarrow+\infty}
u_{kk}^{(n_m)}$ and $\widetilde{u}_{kj}=\displaystyle\lim_{
m\rightarrow+\infty} u_{kj}^{(n_m)}$. Then
$\displaystyle\lim_{m\rightarrow+\infty}U^{(n_{m})}_{i\preceq }=
\widetilde{U}_{i\preceq }$, so that $Z=\widetilde{U}_{i\preceq
}\widetilde{U}_{i\preceq }^{\ast}=(z_{kj})_{k,j \in I}$. As
$z_{ii}\geq 0$, we will consider separately the case $z_{ii}=0$ and
the case $z_{ii}>0$.

\n Suppose that $z_{ii}=0.$ Then
$\widetilde{u}_{ii}=(z_{ii})^{1/2}=0, $ so that
$z_{ki}=\widetilde{u}_{ii}\widetilde{u}_{ki}=0, \ i\prec k$. Thus
$Z=Z_{i\prec}\in \overline{{\cal{P}}}_{i \prec}$, and the result
follows according to (\ref{p_iprec}).

If  $z_{ii}>0$, we consider the elements of  ${\cal{A}}_{i \preceq}$
$${\widetilde{u}}^{i}=\left(
  \begin{array}{cc}
   \widetilde{u}_{ii} & 0 \\
    \widetilde{C}_{i} & e_{\check{i}} \\
  \end{array}
\right)\ \mbox{and} \ {\widetilde{u}}_{i}=\left(
                             \begin{array}{cc}
                               1 & 0 \\
                               0 & \widetilde{U}_{i\prec} \\
                             \end{array}
                           \right),$$

\n where ${\widetilde{C}}_{i}=\displaystyle\sum_{i\prec
j}\widetilde{u}_{ji}$ and $\widetilde{U}_{i
\prec}=(\widetilde{u}_{jk})_{j, k \in I_{i\prec}}$. Then
$\widetilde{u}_{i\preceq }={\widetilde{u}}^{i}{\widetilde{u}}_{i}$
and $\widetilde{u}_{ii}=(z_{ii})^{1/2}>0$. Let $T_1$ and $T_2$ in
${\cal{T}}_{l}^{+}$, such that $T_{1}={\widetilde{u}}^{i} $ in
${\cal{A}}_{i \preceq}$ and $T_2={\widetilde{u}}_{i}$ in
${\cal{A}}_{i \preceq}$, we have $T_{2}.\check{E}^{\check{i}}\in
\overline{{\cal{P}}}_{i\prec}$. By induction hypothesis, there
exists a unique $\psi_1\in \varepsilon^{i}$ such that $\psi_1(i)=0$
and there exists $\tilde{T}_2 \in {\cal{T}}_{l}^{+}$ such that
$T_2.\check{E}^{\check{i}}=\tilde{T}_2.e_{\psi_1}$. Let $\psi_{2}\in
\varepsilon^{i}$, such that $\psi_2(i)=1$ and $\psi_2(j)=0 \ \forall
j \not=i$ and put $\psi =\psi_{1}+\psi_{2}\in \varepsilon^{i}$ and
$T=T_1\tilde{T}_2\in {\cal{T}}_{l}^{+}$. Then we have
\begin{eqnarray*}
Z&=&T_1T_2.(E_i+\check{E}^{\check{i}})\nonumber\\
&=&T_1.E_i+T_{2}.\check{E}^{\check{i}}\nonumber\\
&=&T_1.E_i+\tilde{T}_{2}.e_{\psi_1}\nonumber\\
&=&(T_1\tilde{T}_2).e_{\psi}\\
&=&T.e_{\psi},
\end{eqnarray*}
and (\ref{ipns p_i bar}) is proved.
\end{Pff}
\subsection{Gamma functions}
We use the generalized power function to introduce a generalized
gamma function on an homogeneous cone.

For $i \in \wp\cup \cal{S}$, $\psi \in \varepsilon^{i}$ and
$\chi_i=\{\lambda_j, \ j \in I; \ \lambda_j=0, \forall j \not\in
I_{i\preceq}\}$, we set
\begin{equation}{\mathcal{X}}(\psi)= \{\chi_i\in {\cal{X}}\ |\  \lambda_{j}=0, \ \textrm{for all} \ j
\in I_{i\preceq} \ \textrm{such that} \ \psi(j)=0\}.\end{equation}
For every $\chi_i \in {\mathcal{X}}(\psi)$, we define a generalized
power function on ${\mathcal{T}}_{l}^{+}.e_{\psi}$ by
\begin{equation}\label{gener power}
\Delta_{\chi_i}^{\psi}(T.e_{\psi})=\Delta_{\chi_i}^{(i)}(TT^{\ast}),
\hskip1.2cm \forall T \in {\mathcal{T}}_{l}^{+},
\end{equation}
where $\Delta_{\chi_i}^{(i)}(TT^{\ast})$ is defined by (\ref{genra
power p_i}). We also define $n^{\psi}=(n^{i}_{j.}, j \in I)$ by

\begin{eqnarray}n^{i}_{j.}=\displaystyle\sum_{k\prec
j}\psi(k)n_{kj}\hskip1.5cm \forall j \in I_{i\preceq}.
\end{eqnarray}

\n When $\psi\not=0$, we introduce the measure $\nu_{\psi}$ on
${\mathcal{T}}_{l}^{+}.e_{\psi}$ defined by
\begin{eqnarray}\label{def measure}
\nu_{\psi}(d(T.e_{\psi}))=\Delta_{{\dot{\chi}_i}^{\psi}}^{(i)}(TT^{\ast})\displaystyle\prod_{\begin{array}{c}
                                                                         i\preceq j \preceq k\\
                                                                         \psi(j)=1

                                                                        \end{array}
}dt_{kj},
\end{eqnarray}

\n where ${\dot{\chi}}_i^{\psi}=\{\lambda_j \in \reel, j \in I, \
\textrm{such that} \ \lambda_j=-\psi(j)(1+n^{i}_{j.})/2, \
\textrm{if} \ j \in I_{i\preceq}\ \textrm{and} \ \lambda_j=0 \
\textrm{if} \ j \not\in I_{i\preceq}\}$, and $T=(t_{jk})_{j,k \in I}
\in {\mathcal{T}}_{l}^{+}$. Finally, we denote by $\nu_{0}$ be the
Dirac measure at $0$.
\begin{theorem}\label{gamma function} Let $i \in \wp\cup
{\cal{S}}$ and $\chi_i=\{\lambda_j, \ j \in I; \ \lambda_j=0,
\forall j \not\in I_{i\preceq}\}\in {\mathcal{X}}(\psi)$. The
integral
\begin{eqnarray}\label{integral gamma cte}
\Gamma_{{\mathcal{T}}_{l}^{+}.e_{\psi}}(\chi_i)=\displaystyle\int_{{\mathcal{T}}_{l}^{+}.e_{\psi}}\exp\{-\mbox{tr}(Z)\}\Delta_{\chi_{i}}^{\psi}(Z)\nu_{\psi}(dZ)
\end{eqnarray}
converges if and only if $\chi_i\in {\mathcal{X}}(\psi)$ satisfies
the following condition:
\begin{eqnarray}\label{condition gamma cte}
\lambda_j > \frac{n^{i}_{j.}}{2} \hskip1.5cm \forall j \in I \   \
\textrm{such that}\   \ \psi(j)=1.
\end{eqnarray}
Moreover, under this condition, one has
\begin{equation}\label{formule gamma cte}
\Gamma_{{\mathcal{T}}_{l}^{+}.e_{\psi}}(\chi_i)= 2^{-|\psi|}
\pi^{-|n^{\psi}|/2}\displaystyle\prod _{\begin{array}{c}
                                                    j \in I \\
                                                    \psi(j)=1
                                                  \end{array}}
\Gamma(\lambda_j-\frac{n^{i}_{j.}}{2}),
\end{equation}
where $|\psi|=\displaystyle\sum_{j \in I_{i\preceq}}
 \psi(j)$ and $|n^{\psi}|=\displaystyle\sum_{j \in
I_{i\preceq}}n^{i}_{j.}\cdot$
\end{theorem}

\begin{Pff} \ If $\psi=0$, the integral (\ref{integral gamma cte}) reduces to
$1$. Thus (\ref{condition gamma cte}) and (\ref{formule gamma cte})
hold trivially. If $\psi\not=0$, then writing $Z=U.e_{\psi}$, where
$U=(u_{jk})_{j,k \in I}\in {\mathcal{T}}_{l}^{+}$, the integral
(\ref{integral gamma cte}) can be written
\begin{eqnarray*}
\Gamma_{{\mathcal{T}}_{l}^{+}.e_{\psi}}(\chi_i
)=\displaystyle\int_{{\mathcal{T}}_{l}^{+}e_{\psi}}\exp\displaystyle\{-\displaystyle(\displaystyle\sum_{\psi(j)=1}
(u_{jj}^{2}+\displaystyle\sum_{j\prec k  }\parallel
u_{kj}\parallel_{kj}^{2}))\}\displaystyle\prod_{\begin{array}{c}
                                                  i\preceq j\prec k   \\
                                                  \psi(j)=1
                                                \end{array}}
u_{jj}^{2\lambda_{j}-n_{j.}^{i}-1}du_{jj}du_{kj}\cdot
\end{eqnarray*}
For $j \in I_{i\preceq}$, let \begin{eqnarray} \label{C_i}
\hskip-0.75cm C_j=(s_{ln})_{l,n \in I}\in
{\cal{C}}^{j}=\displaystyle\sum_{j\prec k}E_{kj}, \ \mbox{with} \
s_{lj}=u_{lj} \ \mbox{if} \ j\prec l \ \mbox{and} \ s_{ln}=0 \
\mbox{otherwise}.
\end{eqnarray}
It is clear that $\parallel
C_j\parallel^{2}=\displaystyle\sum_{j\prec k}\parallel u_{kj}
\parallel ^{2}_{kj}\ \textrm{and} \ dC_j=\displaystyle\prod_{j\prec
k}du_{kj}.$ Hence

$$\Gamma_{{\mathcal{T}}_{l}^{+}.e_{\psi}}(\chi_i)=
\displaystyle\prod_{{\begin{array}{c}
                             i\preceq j\\
                             \psi(j)=1
                           \end{array}}}\displaystyle\int_{0}^{+\infty}\exp^{-u_{jj}^{2}}u_{jj}^{2\lambda_{j}-n_{j.}^{i}-1}du_{jj}\displaystyle\prod_{{\begin{array}{c}
                             i\preceq j\\
                             \psi(j)=1
                           \end{array}}}\displaystyle\int_{{\cal{C}}^{j}}\exp^{\{-\parallel C_j\parallel
                           ^{2}\}}dC_j.$$

\n Therefore the convergence condition is reduced to the one
corresponding to the ordinary gamma functions, that is $\lambda_j
> \frac{n^{i}_{j.}}{2} \ \forall j \in I \ \textrm{such that}\   \
\psi(j)=1.$
\end{Pff}
\begin{remark}
From Theorem \ref{gamma function}, we have, for $i \in I$, a
relation between $\Gamma_{{\mathcal{T}}_{l}^{+}.e_{\psi}}(\chi_{i})$
and $\Gamma_{\cal{P}}({\chi})$. In fact, if we denote by
$\textbf{1}_{i}$ the element of $\varepsilon^{i}$, such that
$\psi(j)=1, \ \forall \ j \in I_{i\preceq}\setminus\mathcal{S},$  if
$\ i \in {\wp}$ and such that $\psi(j)=1, \ \forall \ j \in
I_{i\preceq},$ if $i \in\mathcal{S}$, then it is clear that
\begin{eqnarray*}{\cal{P}}=\displaystyle\sum_{i \in \wp\cup
{\cal{S}}}{\mathcal{T}}_{l}^{+}.e_{\bf{1}_{i}},\end{eqnarray*} and
using (\ref{les n_i}), for $\chi=\displaystyle\sum_{i \in
I}\chi_i\in{\cal{X}}$, where $\chi_i=\{\lambda_j, \ j \in I; \
\lambda_j=0, \forall j \not\in I_{i\preceq}\}$,  we have
\begin{eqnarray*}\displaystyle\prod_{i \in
\wp\cup
\mathcal{S}}\Gamma_{{\mathcal{T}}_{l}^{+}.e_{\bf{1}_{i}}}(\chi_{i})=2^{-|I|}\pi^{\frac{n_{.}-|I|}{2}}\displaystyle\prod_{i
\in I}
\Gamma(\lambda_i-\frac{n_{i.}}{2})=2^{-|I|}\Gamma_{\cal{P}}({\chi})\cdot
\end{eqnarray*}
\end{remark}
\subsection{Riesz measures}
For the definition of the Riesz distribution, we need to introduce
some other notations. Let $i \in \wp\cup \cal{S}$ and
$\chi_i=\{\lambda_j, \ j \in I; \ \lambda_j=0, \forall j \not\in
I_{i\preceq}\}$ and introduce for $\psi \in \varepsilon^{i}$ and
$\omega\in \varepsilon^{\check{i}}$, the following sets
\begin{eqnarray}\label{B(i,psi)}\hskip-0.5cm{\cal{B}}(i,\psi)=\displaystyle\left\{\chi_{i} \in {\cal{X}};\ \lambda_{j}=0 \
\textrm{for} \ j \not\in I_{i\preceq},
\lambda_{j}=\frac{n^{i}_{j.}}{2} \ \textrm{when}
 \ j \in I_{i\preceq} \ \textrm{and} \ \psi(j)=0  \right\},\end{eqnarray}
\begin{eqnarray}\label{check
B(i,psi)}\hskip-0.5cm{\check{\cal{B}}}(\check{i},\omega)=\left\{\chi_i\in
{\cal{X}};\ \upsilon_{j}=0 \ \textrm{for} \ j \not\in I_{i\prec},
\vartheta_{j}=\frac{n^{i}_{j.}}{2} \ \textrm{when}
 \ j \in I_{i\prec} \ \textrm{and} \ \omega(j)=0 \right\},\end{eqnarray}
\begin{eqnarray}\label{{Xi}(i,psi)}
{\Xi}(i, \ \psi)=\displaystyle\left\{\chi_{i} \ \in {\cal{B}}(i, \
\psi), \ \lambda_{j}>\frac{n^{i}_{j.}}{2}\ \textrm{when} \  j \in
I_{i\preceq}  \ \textrm{and} \ \psi(j)=1 \right\},
\end{eqnarray}
\begin{eqnarray}\label{{check{Xi}}(i,psi)}{\check{\Xi}}(\check{i}, \ \omega)=\left\{\chi_{i} \ {\check{\cal{B}}}(\check{i},\omega)
\ \lambda_{j}>\frac{n^{i}_{j.}}{2}\ \textrm{when} \  j \in
I_{i\prec} \ \textrm{and} \ \omega(j)=1\right\},\end{eqnarray}
\begin{eqnarray} \label{Xi(i)}{\Xi}(i)=\bigcup_{\psi \in
\varepsilon^{i}}{\Xi}(i, \ \psi), \hskip1cm
\check{\Xi}(\check{i})=\displaystyle\bigcup_{\omega \in
\varepsilon^{\check{i}}}{\check{\Xi}}(\check{i}, \
\omega),\hskip1cm\Xi=\displaystyle\sum_{i \in \wp\cup
\cal{S}}{\Xi}(i)\cdot
\end{eqnarray}
\n For every $\chi_{i}=\{\lambda_j, \ j \in I\}\in {\cal{X}}$, let
\begin{eqnarray}\label{chi tide}
\widetilde{\chi}_{i}=\{\lambda_j-(1-\psi(j))\frac{n_{j.}^{i}}{2},
\textrm{when}\ j \in I_{i\preceq}, \ \textrm{and} \ 0 \ \textrm{if}
\  j \not\in I_{i\preceq}\}.
\end{eqnarray}
It is clear that if $\chi_{i} \in {\cal{B}}(i, \ \psi)$, then
$\widetilde{\chi}_{i} \in {\mathcal{X}}(\psi)$.

In what follows, we denote the Laplace transform of a positive
measure $\mu$ on the cone ${\cal {P}}$ by
\begin{eqnarray}
L_{\mu}(\theta)=\displaystyle\int_{\cal{P}}\exp\{-\textrm{tr}(\theta
Z)\}\mu(dZ),  \hskip1cm \theta\in {\cal{P}}^{\ast}.
\end{eqnarray}
\begin{theorem}\label{main result} There exists a positive measure $R_{\chi}$ such
that the Laplace transform is defined on ${\cal{P}}^{\ast}$ and is
equal to $\Delta_{\chi}(\theta^{-1})$ if and only if $\chi\in \Xi$.
\end{theorem}
The proof  of Theorem \ref{main result} relies on the following
proposition.
\begin{prop} \label{prop main result} Let $i \in \wp\cup \mathcal{S} $. Then there exists a positive measure $R_{\chi_i}$ such that the Laplace
transform is defined on ${\cal{P}}^{\ast}$ and is equal to
$\Delta^{(i)}_{\chi_{i}}(\theta^{-1})$ if and only if $\chi_{i}\in
\Xi(i)$.
\end{prop}

\begin{Pff}\ $\Leftarrow)$ Let $\chi_{i}\in
\Xi(i)$. Then there exists $\psi \in \varepsilon^{i}$ such that
${\chi}_{i} \in {\Xi}(i, \ \psi)$. It is clear that
$\widetilde{\chi}_{i}$ defined by (\ref{chi tide}) satisfies
(\ref{condition gamma cte}). We will show that the Laplace transform
of the measure
\begin{eqnarray*}
R_{\chi_{i}}(dZ)=\frac{1}{\Gamma_{{\mathcal{T}}_{l}^{+}.e_{\psi}}(\widetilde{\chi}_{i})}
\displaystyle\Delta_{{\widetilde{\chi}}_{i}}^{\psi}(Z){\bf
1}_{{\mathcal{T}}_{l}^{+}.e_{\psi}}(Z)\nu_{\psi}(dZ)
\end{eqnarray*}
is defined on ${\cal{P}}^{\ast}$ and is given by
\begin{eqnarray*}L_{R_{\chi_{i}}}(\theta)=\frac{1}{\Gamma_{{\mathcal{T}}_{l}^{+}.e_{\psi}}(\widetilde{\chi}_{i})}
\displaystyle\int_{{\mathcal{T}}_{l}^{+}.e_{\psi}}\exp\{-\mbox{tr}(\theta
Z)\}\Delta_{{\widetilde{\chi}}_{i}}^{\psi}(Z)\nu_{\psi}(dZ)=\Delta^{(i)}_{\chi_{i}}(\theta^{-1}).\end{eqnarray*}
In fact, as $\theta \in {\cal{P}}^{\ast}$, then $\theta_{i\preceq}$
defined by (\ref{Z_iprece}) is in the dual cone ${\cal{P}}_{i
\preceq}^{\ast}$ of ${\cal P}_{i\preceq}$, and there exists $T$
 in ${\mathcal{T}}_{l}^{+}$ such that
$\theta_{i\preceq}=T^{\ast}.E^{i}$. Let $Y=\pi(T)(Z)$, where $\pi$
is defined by (\ref{pi}). As $Z \in {\mathcal{T}}_{l}^{+}.e_{\psi}$,
there exists $S \in {\mathcal{T}}_{l}^{+}$ such that $Z=S.e_{\psi}$.
This with (\ref{def measure}) imply that
\begin{eqnarray}\label{dnu(y)}
\nu_{\psi}(dY)=\nu_{\psi}(d(\pi(T)(Z)))=\Delta^{(i)}_{{\ddot{\chi}_i}^{\psi}}(\theta)
\nu_{\psi}(dZ),\end{eqnarray} where ${\ddot{\chi}}_{i}^{\psi
}=\{\lambda_j \in \reel, j \in I \ \textrm{such that} \
\lambda_{j}=(1-\psi(j))\frac{n^{i}_{j.}}{2}, \ \textrm{if} \ j \in
I_{i\preceq}\ \textrm{and} \ 0 \ \textrm{if} \ j \not\in
I_{i\preceq}\}$.

\n Since $\widetilde{\chi}_i \in {\mathcal{X}}(\psi)$, then
\begin{eqnarray}\label{Delta(x)}
\Delta_{{\widetilde{\chi}_i}}^{\psi}(Z)&=&\Delta_{{\widetilde{\chi}_i}}^{\psi}(\pi^{-1}(T)(Y))\nonumber\\
&=&\Delta_{{\widetilde{\chi}_i}}^{\psi}(T^{-1}.e_{\psi})\Delta_{{\widetilde{\chi}_i}}^{\psi}(Y)\nonumber\\
&=&\Delta^{(i)}_{\widetilde{\chi}_i}(\theta^{-1})\Delta_{{\widetilde{\chi}_i}}^{\psi}(Y).
\end{eqnarray}
Using (\ref{chi tide}), (\ref{dnu(y)}) and (\ref{Delta(x)}), we get
\begin{eqnarray}
\Delta_{{\widetilde{\chi}_i}}^{\psi}(Y)\nu_{\psi}(dY)=\Delta^{(i)}_{\chi_{i}
}(\theta)\Delta_{{\widetilde{\chi}_i}}^{\psi}(Z)\nu_{\psi}(dZ).
\end{eqnarray}
Then
\begin{eqnarray*}L_{R_{\chi_{i}}}(\theta)&=&\frac{1}{\Gamma_{{\mathcal{T}}_{l}^{+}.e_{\psi}}(\widetilde{\chi}_{i})}
\displaystyle\int_{{\mathcal{T}}_{l}^{+}.e_{\psi}}\exp\{-\mbox{tr}(\theta
\pi^{-1}(T)(Y))\}\Delta^{(i)}_{\chi_i}(\theta^{-1})\Delta_{{\widetilde{\chi}}_{i}}^{\psi}(Y)\nu_{\psi}(dY)\\
&=&\Delta^{(i)}_{\chi_i}(\theta^{-1})\frac{1}{\Gamma_{{\mathcal{T}}_{l}^{+}.e_{\psi}}(\widetilde{\chi}_{i})}
\displaystyle\int_{{\mathcal{T}}_{l}^{+}.e_{\psi}}\exp\{-\mbox{tr}
Y\}\Delta_{{\widetilde{\chi}}_{i}}^{\psi}(Y)\nu_{\psi}(dY)\\
&=&\Delta^{(i)}_{\chi_i}(\theta^{-1}).
\end{eqnarray*}
\n $\Rightarrow)$ Suppose that there exists a positive measure
$R_{\chi_{i}}$ such that the Laplace transform is defined on
${\cal{P}}^{\ast}$ and is equal to
$\Delta^{(i)}_{\chi_i}(\theta^{-1})$.  Our aim to show that
$\chi_{i} \in {\Xi}(i)$.

\n For this $\chi_i$ and a $\psi$ in $\varepsilon^{i}$, consider the
generalized positive Riesz measure which we also denote
$R_{\chi_{i}}$ defined for $\varphi$ in the Schwartz space
$\mathcal{S}(\cal{A})$ of rapidly decreasing functions on $\cal{A}$
by
\begin{eqnarray}\label{integral riesz}
R_{\chi_{i}}(\varphi)=\frac{1}{\Gamma_{{\mathcal{T}}_{l}^{+}.e_{\psi}}(\widetilde{\chi}_{i})}
\displaystyle\int_{{\mathcal{T}}_{l}^{+}.e_{\psi}}\varphi(Z)\Delta_{{\widetilde{\chi}}_{i}}^{\psi}(Z)\nu_{\psi}(dZ).
\end{eqnarray}
We will prove by induction on the rank of the cone ${\cal{P}}_{i
\preceq}$ that $\chi_{i} \in {\Xi}(i, \ \psi)$. Suppose that
$\mbox{rank}{\cal{P}}_{i\preceq}=1.$ Then we have either cardinality
of $\wp$ equal to 1 or cardinality of $\mathcal{S}$ equal to 1. Thus
$R_{\chi_{i}}$ coincides with the Riesz measure $\rho_{\lambda}$ on
$]0,+\infty[$ given by
\begin{eqnarray}\label{rho_lambda}{\rho_{\lambda}}(\varphi)=\frac{1}{\Gamma(\lambda)}\displaystyle\int_{0}^{+\infty}\varphi(u)u^{\lambda-1}du.\end{eqnarray}
This implies that $\chi_i\equiv\lambda$ and $\lambda>0$ which means
that the result is true when $\mbox{rank}{\cal{P}}_{i\preceq}=1$.
Now, suppose that the claim holds for any $i\in \wp\cup \cal{S}$
such that $\mbox{rank}{\cal{P}}_{i\preceq}\leq k-1$, and let us show
that it also holds  for $i$ such that
$\mbox{rank}{\cal{P}}_{i\preceq}=k$. Consider $i \in \wp\cup
\cal{S}$ such that $\mbox{rank}{\cal{P}}_{i\preceq}=k$. Then
$\mbox{rank}{\cal{P}}_{i\prec}=k-1$ so that
$\mbox{rank}{\cal{P}}_{j\preceq}\leq k-1$, $\forall j \in M_i$.
Using the induction hypothesis, we have that $\forall j \in M_i$,
$\xi_{j}=\{\beta_l\in \reel, \ l \in I, \ \beta_l=0; \ \forall l
\not \in I_{j\preceq} \}$ is in  $\Xi(j)$. Let ${R}_{{\xi}_{j}}$ be
a Riesz measure defined as in (\ref{integral riesz}) on the cone
${\cal{P}}_{j \preceq}$ for some $\psi$ in $\varepsilon^{j}$ and let
$\check{\xi}_{\check{i}}=\displaystyle\sum_{j \in
M_i}\xi_{j}=\{\beta_j \in \reel, j \in I; \beta_{j}=0 \ \textrm{for}
\ j \not\in I_{i\prec}\}$. Then from (\ref{P_M}), the measure
$\check{R}_{\check{\xi}_{\check{i}}}=\displaystyle\prod^{\ast}_{j
\in M_i}R_{{\xi}_{j}}$, where $\displaystyle\prod^{\ast}$ is the
convolution product, is concentrated on ${\cal{P}}_{i \prec}$.
Consider the sets
$$\check{\chi}_{\check{i}}=\{\lambda_j \in \reel, j \in I;
\lambda_{j}=0 \ \textrm{for} \ j \not\in I_{i\prec}\},$$
$$\check{n}^{i}=\{\beta_{k}\in \reel, \ k \in I; \beta_{k}=n_{ki}\
\textrm{for} \ k \in I_{i\prec} \  \textrm{and} \ \beta_{k}=0 \
\textrm{for} \ k \not\in I_{i\prec}\},$$ and
$$M(\lambda_{i})=\{\alpha_{k}, \ k \in I\},$$ where
$$\alpha_{k}=\left \{\begin{array}{ccccc}\lambda_i & \textrm{for} & k=i&&\\
\frac{n_{ki}}{2} & \textrm{for}&k \in I_{i\prec}&& \\
 0 &  & \textrm{otherwise.}&& \end{array}\right.$$
 Then it is easy to verify that $M(\lambda_{i})\in {\cal{B}}(i,\psi_{1})$, where $\psi_{1}
\in \varepsilon^{i}$ such that $\psi_{1}(i)=1$, and $\psi_{1}(j)=0 $
$\forall j\not=i$. Also we
 have
$\chi_{i}-M(\lambda_{i})=\check{\chi}_{\check{i}}-\frac{\check{n}^{i}}{2}\in
\check{{\cal{B}}}(\check{i},\omega)$, where
$\check{{\cal{B}}}(\check{i},\omega)$ is defined by (\ref{check
B(i,psi)}).

\n Using the Laplace transforms, we obtain that
\begin{eqnarray}\label{convolution}R_{\chi_{i}}=R_{M(\lambda_{i})}\ast
\check{R}_{\check{\chi}_{\check{i}}-\frac{\check{n}^{i}}{2}}.\end{eqnarray}
 Proceeding as in the proof
of Theorem \ref{gamma function}, and using (\ref{formule gamma
cte}), we get
$$R_{M(\lambda_{i})}(\varphi)=\frac{2
\pi^{-\frac{\dim{\cal{C}}^{i}}{2}}}{\Gamma(\lambda_i)}
\displaystyle\int_{{\mathcal{T}}_{l}^{+}e_{\psi_1}}\varphi(U.e_{\psi_1})u_{ii}^{2\lambda_{i}-1}du_{ii}dC_i,$$
where $U\in{\mathcal{T}}_{l}^{+}$, $U.e_{\psi_1}=(u_{jk})_{j,k \in
I}$ and $C_i$ is defined by (\ref{C_i}). On the other hand, it is
easy to verify that
\begin{eqnarray*}U.e_{\psi_1}=u_{ii}^{2}E_i+u_{ii}C_i+\displaystyle\sum_{i\prec j}\parallel
u_{ji}\parallel_{ji}^{2}E_j.
\end{eqnarray*}
Hence, if we define
\begin{eqnarray*}
{\mathcal{Q}}^{i}\;:\;{{\cal{C}}}^{i}\times{{\cal{C}}}^{i}&\rightarrow&{\cal{A}}_{i\prec}\nonumber\\
(C_i,C_i)&\mapsto&{\mathcal{Q}}^{i}(C_i,C_i)=\displaystyle\sum_{i\prec
j}\parallel u_{ji}\parallel_{ji}^{2}E_j,
\end{eqnarray*} then we have
\begin{eqnarray*}U.e_{\psi_1}=u_{ii}^{2}E_i+u_{ii}C_i+{\mathcal{Q}}^{i}(C_i,C_i).\end{eqnarray*}
Setting $u_{ii}=\sqrt{v}$, we get
$$R_{M(\lambda_{i})}(\varphi)=\frac{2
\pi^{-\frac{\dim{\cal{C}}^{i}}{2}}}{\Gamma(\lambda_i)}\displaystyle\int_{0}^{+\infty}
\displaystyle\int_{{{\cal{C}}}^{i}}\varphi(vE_i+\sqrt{v}C_i+{\mathcal{Q}}^{i}(C_i,C_i))dC_iv^{\lambda_{i}-1}dv.$$
This, using (\ref{rho_lambda}), becomes
\begin{eqnarray}\label{R_M}R_{M(\lambda_{i})}(\varphi)= \pi^{-\frac{\dim{\cal{C}}^{i}}{2}}\rho_{\lambda_{i}}(
\displaystyle\int_{{{\cal{C}}}^{i}}\varphi(vE_i+\sqrt{v}C_i+{\mathcal{Q}}^{i}(C_i,C_i))dC_i)_{v}.
\end{eqnarray}
As for $\lambda_i=0$, $\rho_{0}$ is the Dirac measure at $v=0$, we
get
\begin{eqnarray}\label{R_M0}R_{M(0)}(\varphi)=R_{\frac{\check{n}^{i}}{2}}(\varphi)=
\pi^{-\frac{\dim{\cal{C}}^{i}}{2}}\displaystyle\int_{{{\cal{C}}}^{i}}\varphi({\mathcal{Q}}^{i}(C_i,C_i))dC_i.
\end{eqnarray}
Using (\ref{convolution}) and (\ref{R_M}), we obtain
\begin{eqnarray}\label{R_chi(phi)}
R_{\chi_{i}}(\varphi)=\pi^{-\frac{\dim{\cal{C}}^{i}}{2}}\rho_{\lambda_{i}}(\displaystyle\int_{{{\cal{C}}}^{i}}\check{R}_{\check{\chi}_i-\frac{\check{n}^{i}}{2}}(\varphi(vE_i+\sqrt{v}C_i+{\mathcal{Q}}^{i}(C_i,C_i)+Y))_{Y}dC_i))_{v},
\end{eqnarray}

\n Denote by $C_{c}^{\infty}$ the set of $C^{\infty}$ functions with
compact support and consider the functions of the form
$$\varphi(Z)=\varphi_{1}(z_{11})\varphi_{2}(Z_{i\prec})$$
where $Z=(z_{ij})_{i,j \in I} \in {\cal{A}}$, $Z_{i\prec} \in
{\cal{A}}_{i\prec},$ $\varphi_{1}\in C_{c}^{\infty}(\reel)$ and
$\varphi_{2}\in C_{c}^{\infty}({\cal{A}}_{i\prec})$. Then by
(\ref{R_M0}) and (\ref{R_chi(phi)}), we have
\begin{eqnarray}\label{r_chi convolution}
R_{\chi_{i}}(\varphi)&=&\pi^{-\frac{\dim{\cal{C}}^{i}}{2}}\rho_{\lambda_{i}}(\varphi_{1})\displaystyle\int_{{{\cal{C}}}^{i}}\check{R}_{\check{\chi}_i-\frac{\check{n}^{i}}{2}}(\varphi_2({\mathcal{Q}}^{i}(C_i,C_i)+Y))_{Y}dC_i\nonumber\\
&=&\rho_{\lambda_{i}}(\varphi_{1})(\check{R}_{\frac{\check{n}^{i}}{2}}\ast\check{R}_{\check{\chi}_i-\frac{\check{n}^{i}}{2}})(\varphi_{2})\nonumber\\
&=&\rho_{\lambda_{i}}(\varphi_{1})\check{R}_{\check{\chi}_{\check{i}}}(\varphi_{2}).
\end{eqnarray}
For a suitable choice of non-negative $\varphi_1\in
C_{c}^{\infty}(\reel)$, we have $\rho_{\lambda_{i}}(\varphi_{1})>0$.
If $\varphi_{2}\geq 0$, then using (\ref{R_chi(phi)}) and the
positivity of $R_{\chi_{i}}$, we get
$\check{R}_{\check{\chi}_{\check{i}}}(\varphi_{2})=(\rho_{\lambda_{i}}(\varphi_{1}))^{-1}R_{\chi_{i}}(\varphi)\geq0$.
Thus $\check{R}_{\check{\chi}_i}$ is positive and the induction
hypothesis ensures that $\check{\chi}_{\check{i}}\in
{\check{\Xi}}(\check{i}, \ \omega)$.

\n Now, fix a non-negative $\varphi_{2}$ such that
$\check{R}_{\check{\chi}_{\check{i}}}(\varphi_{2})$ is strictly
positive. Then using again (\ref{R_chi(phi)}), we get
$\rho_{\lambda_{i}}(\varphi_{1})\geq0$ for any $\varphi_{1}\geq0$.
Therefore $\rho_{\lambda_{i}}$ is positive and we deduce that
$\lambda_{i}\geq0$. If $\lambda_{i}=0,$ then choosing a $\psi$ in
$\varepsilon^{i}$ such that $\psi(i)=0$, we get $\chi_{i} \in
{\Xi}(i, \ \psi)$. To study the case $\lambda_i>0$, we first observe
that the map
\begin{eqnarray*}
]0,+\infty[\times{{\cal{C}}}^{i}\times{\cal{A}}_{i\prec}&\rightarrow&\{X \in {\cal{A}}, x_{11}>0\}\nonumber\\
(v,C_i,Y)&\mapsto&vE_i+\sqrt{v}C_i+{\mathcal{Q}}^{i}(C_i,C_i)+Y
\end{eqnarray*} is a diffeomorphism whose the inverse is given by $$ x\mapsto(x_{11},\
x_{11}^{-1/2}\displaystyle\sum_{1\prec k}X_{k1}, \
x_{i\prec}-\frac{1}{x_{11}}
{\mathcal{Q}}^{i}(\displaystyle\sum_{1\prec
k}X_{k1},\displaystyle\sum_{1\prec k}X_{k1}).$$ For a functions the
functions $\varphi \in C_{c}^{\infty}({\cal{A}})$ of the form
\begin{eqnarray*}
\varphi(Z)= \left \{ \begin{array}{cc}
\varphi_{1}(v)\varphi_{2}(C_i)\varphi_{3}(Y) & (z_{11}>0),\\
0 & (z_{11}\leq0),
\end{array}\right.
\end{eqnarray*}
with $\varphi_{1}\in  C_{c}^{\infty}(]0,+\infty[)$, $\varphi_{2}\in
C_{c}^{\infty}({{\cal{C}}}^{i})$ and $\varphi_{3}\in
C_{c}^{\infty}({\cal{A}}_{i\prec})$, by (\ref{R_chi(phi)}), we have
that
\begin{eqnarray*}
R_{\chi_{i}}(\varphi)=\pi^{-\frac{\dim{\cal{C}}^{i}}{2}}\rho_{\lambda_{i}}(\varphi_{1})\check{R}_{\check{\chi}_i-\frac{\check{n}^{i}}{2}}(\varphi_{3})\displaystyle\int_{{{\cal{C}}}^{i}}\phi_{2}(C_i)dC_i.\nonumber\\
\end{eqnarray*}
Since $\lambda_{i}>0,$ the positivity assumption of $R_{\chi_i}$
yields that $\check{R}_{\check{\chi}_i-\frac{\check{n}^{i}}{2}}$ is
positive. This by the induction hypothesis implies that
$\check{\chi}_i-\frac{\check{n}^{i}}{2}\in {\check{\Xi}}(\check{i},
\ \omega)$. Finally, choose a $\psi$ in $\varepsilon^{i}$ such that

$$\psi(j)=\left\{\begin{array}{ccc}\omega(j)& \textrm{for}& \forall j\not=i\\
1& \textrm{for}&j=i, \end{array}\right.$$ where $\omega \in
\varepsilon^{\check{i}}$. Then, as
$\check{\chi}_i-\frac{\check{n}^{i}}{2}\in {\check{\Xi}}(\check{i},
\ \omega)$, for $j \in I_{i\prec}$, we have that
$\lambda_{j}=\frac{n^{i}_{j.}}{2}$ if $ \omega(j)=0$ and
$\lambda_{j}>\frac{n^{i}_{j.}}{2}$ if $\omega(j)=1$. As $\lambda_{i}
>0$ and $n^{i}_{i.}=0$, then $\lambda_{i}>\frac{n^{i}_{i.}}{2}$. This
means that for such a $\psi$, we have that $\chi_{i} \in {\Xi}(i, \
\psi)$. Hence $\chi_{i} \in {\Xi}(i)$ and Proposition \ref{prop main
result} is proved.
\end{Pff}

\begin{Pff} {\bf of Theorem \ref{main result}} \  $(\Leftarrow)$ Let
$\chi=\displaystyle\sum_{i \in \wp\cup \mathcal{S}}\chi_{i}\in \Xi$,
where $\chi_{i} \in {\Xi}(i)$ and let
\begin{eqnarray}\label{R_chi}
R_{\chi}=\displaystyle\prod^{\ast}_{i \in \wp\cup
\mathcal{S}}R_{\chi_{i}},
\end{eqnarray}
$R_{\chi_{i}}$ is the positive measure defined from Proposition
\ref{prop main result}. Then, for $\theta \in {\cal{P}}^{\ast}$
\begin{eqnarray*}
L_{R_{\chi}}(\theta)&=&\displaystyle\prod_{i \in \wp\cup
\mathcal{S}}L_{R_{\chi_{i}}}(\theta)\\
&=&\displaystyle\prod_{i \in \wp\cup
\mathcal{S}}\Delta^{(i)}_{\chi_{i}}(\theta^{-1})\\
&=&\Delta^{(i)}_{\small{\displaystyle\sum_{i \in \wp\cup
\mathcal{S}}}\chi_{i}}(\theta^{-1})\\
&=&\Delta_{\chi}(\theta^{-1})\cdot
\end{eqnarray*}
\noindent $(\Rightarrow)$ We have
$L_{R_{\chi}}(\theta)=\Delta_{\chi}(\theta^{-1}),$ then using the
fact $\displaystyle\prod_{i \in \wp\cup
\mathcal{S}}\Delta_{\chi_{i}}^{(i)}(\theta^{-1})
=\Delta_{\chi}(\theta^{-1})$, Proposition \ref{prop main result},
and putting $R_{\chi}=\displaystyle\prod^{\ast}_{i \in \wp\cup
\mathcal{S}}R_{\chi_{i}}$, such that for
$\theta\in{\cal{P}}^{\ast}$,
$L_{R_{\chi_{i}}}(\theta)=\Delta_{\chi_{i}}^{(i)}(\theta^{-1})$, we
get $\chi_{i}\in {\Xi}(i)$. Therefore $\chi=\displaystyle\sum_{i \in
\wp\cup \mathcal{S}}\chi_{i} \in \Xi$.
\end{Pff}

Following the terminology used in the paper by Hassairi and Lajmi
(2001) in the case of symmetric matrices, we call the measures
$R_{\chi}$, defined above in terms of their Laplace transforms,
Riesz measures on the homogeneous cone. These measures are divided
into two classes according to the position of $\chi \in \Xi$. A
class of measures which are absolutely continuous with respect to
the Lebesgue measure on ${\cal{P}}$ and a class concentrated on
the boundary $\partial{\cal{P}}$ of ${\cal{P}}$.

\begin{prop} Let $\chi=\{\lambda_i, \ i \in I\}\in \cal{X}$.
Then $R_{\chi}$ is absolutely continuous if and only if
$\lambda_{i}>\frac{n_{i.}}{2}, \ i \in I$. In this case
\begin{eqnarray}
\label{wishart}
R_\chi(dZ)=\frac{1}{\Gamma_{\cal{P}}({\chi})}\Delta_{\chi+\ddot{\chi}}(Z){
1}_{\mathcal{P}}(Z)dZ,
\end{eqnarray}
\n where $\ddot{\chi}=\{-n_{i}, \ i \in I\}$ and
$\Gamma_{{\cal{P}}}({\chi})=\pi^{\frac{n_{.}-|I|}{2}}\prod
\limits_{i\in I}\Gamma(\lambda_{i}-\frac{n_{i.}}{2})$.
\end{prop}

\begin{Pff}
 \noindent $\Rightarrow)$ We have
${\cal{P}}=\displaystyle\sum_{i \in \wp\cup
\cal{S}}{\mathcal{T}}_{l}^{+}.e_{\bf{1}_{i}}$, then
$\widetilde{\chi}_i$ defined by (\ref{chi tide}) is equal to
$\chi_{i}$.

\n From (\ref{R_chi}), we have
$R_{\chi}=\displaystyle\prod^{\ast}_{i \in \wp\cup
\mathcal{S}}R_{\chi_{i}}$. Writing  $Z=UU^{\ast} \in {\cal{P}}$,
$U\in {\mathcal{T}}_{l}^{+}$, then $Z=\displaystyle\sum_{i \in
I}Z_i$ where $Z_i$ is in ${\mathcal{T}}_{l}^{+}.e_{\bf{1}_{i}}$,
 and using the proof of Proposition
\ref{prop main result}, we have that for $Z_i =U^{1}_{i\preceq
}U^{1\ast}_{i\preceq }$, $U^{1}_{i\preceq }=(u_{kj})_{k,j \in I}$
\begin{eqnarray*}
R_{\chi_{i}}(dZ_i)=\frac{1}{\Gamma_{{\mathcal{T}}_{l}^{+}.e_{\bf{1}_{i}}}(\chi_{i})}
\Delta_{\chi_{i}}(Z_i)\nu(dZ_i),
\end{eqnarray*}
where $\Delta^{\bf{1}_{i}}_{\chi_{i}}=\Delta_{\chi_i}$ and
$\nu(dZ_i)=\nu_{\bf{1}_{i}}(dZ_i)=\Delta^{(i)}_{{\dot{\chi}_i}^{\bf{1}_{i}}}(Z
)\displaystyle\prod_{\begin{array}{c}
                                                                         i \preceq j \preceq k\\
                                                                         \psi(j)=1

                                                                        \end{array}
}du_{kj}.$

\n Then $\nu(dZ)=\Delta_{{\dot{\chi}}}({Z})dU$, where
$${\dot{\chi}}=\{-(\frac{1+n_{j\cdot}}{2}), \ j \in I\} \
\textrm{and} \ dU=\displaystyle\prod_{i \in
\wp\cup\mathcal{S}}\displaystyle\prod_{\begin{array}{c}
                                                                         i \preceq j \preceq k\\
\psi(j)=1\end{array} }du_{kj}.$$ Using the fact that the mapping $
U\in {\cal{T}}_l^+\mapsto UU^\ast\in{\cal{P}}$ is a diffeomorphism,
we have that $dU=2^{-\mid I\mid}\Delta_{\hat{\chi}}(Z)dZ$, where
$\hat{\chi}=\{-\frac{1+n_{\cdot j}}{2}, \ \forall j \in I\}$. As
$n_j=1+\frac{1}{2}(n_{\cdot j}+n_{j\cdot})$, we have
\begin{eqnarray*}
R_{\chi}(dZ)=\frac{1}{\Gamma_{{\cal{P}}}(\chi)}
\Delta_{\chi}(Z)\Delta_{\dot{\chi}}(Z)\Delta_{\hat{\chi}}(Z)dZ=\frac{1}{\Gamma_{{\cal{P}}}(\chi)}
\Delta_{\chi+\ddot{\chi}}(Z)dZ,
\end{eqnarray*}
where $\ddot{\chi}=\{-n_j, \ j \in I\}$ and
${\Gamma_{{\cal{P}}}(\chi)}=2^{\mid I\mid}\displaystyle\prod_{i \in
\wp\cup\mathcal{S}}\Gamma_{{\mathcal{T}}_{l}^{+}.e_{\psi}}(\chi_{i})=\pi^{\frac{n_{.}-|I|}{2}}\displaystyle\prod_{i
\in I} \Gamma(\lambda_i-\frac{n_{i.}}{2})\cdot$ Moreover, the
condition $\lambda_{i}>\frac{n_{i.}}{2}, \ i \in I$ is easily
deduced from Theorem \ref{gamma function}.

\n $\Leftarrow)$ It suffices to verify that for $\chi$ such that
$\lambda_{i}>\frac{n_{i.}}{2}, \ i \in I$, the Laplace transform of
the measure
$$\frac{1}{\Gamma_{\cal{P}}({\chi})}\Delta_{\chi+\ddot{\chi}}(Z){
1}_{\mathcal{P}}(Z)dZ$$ is equal to $\Delta_{\chi}(\theta^{-1})$,
$\forall \theta \in {\cal{P}}^{\ast}$. In fact, let $\theta \in
{\cal{P}}^{\ast}$, then there exists $T=(t_{ij})$ in
${\mathcal{T}}_{l}^{+}$ such that $\theta=T^{\ast}T$.
 \n Let
$Y=\pi(T)(Z)$, where $\pi$ is defined by (\ref{pi}). Then $dZ=\det
\pi^{-1}(T)dY$ and
$$\Delta_{\chi+\ddot{\chi}}(Z)=\Delta_{\chi+\ddot{\chi}}
(\pi^{-1}(T)Y)=\Delta_{\chi+\ddot{\chi}}
(\theta^{-1}Y)=\Delta_{\chi+\ddot{\chi}}(\theta^{-1})\Delta_{\chi+\ddot{\chi}}(Y).$$
 From Andersson and Wojnar (2004), we have $$\det
\pi^{-1}(T)=\displaystyle\prod_{i \in
I}t_{ii}^{-2n_i}=\Delta_{-\ddot{\chi}}(\theta^{-1}),$$ where
$-\ddot{\chi}=\{n_i, i \in I\}$. Writing $Y=\displaystyle\sum_{i \in
I}Y_i$, where $Y_i$ is in ${\mathcal{T}}_{l}^{+}.e_{\bf{1}_{i}}$ ,
then as
\begin{eqnarray*}\Gamma_{\cal{P}}({\chi})=2^{|I|}\displaystyle\prod_{i \in
\wp\cup
\mathcal{S}}\Gamma_{{\mathcal{T}}_{l}^{+}.e_{\bf{1}_{i}}}(\chi_{i})&=&2^{|I|}\displaystyle\prod_{i
\in \wp\cup
\mathcal{S}}\displaystyle\int_{{\mathcal{T}}_{l}^{+}.e_{\bf{1}_{i}}}\exp\{-\mbox{tr}Y_i\}\Delta_{\chi_{i}}^{\bf{1}_{i}}(Y_i)\nu_{\bf{1}_{i}}(dY_i)\\
&=&\int_{\cal{P}}\exp-\{\textrm{tr}Y\}\Delta_{\chi+\ddot{\chi}}(Y)dY,
\end{eqnarray*}
we obtain that \begin{eqnarray*}\frac{1}{\Gamma_{{\cal{P}}}(\chi)}
\int_{\cal{P}}\exp\{-\textrm{tr}(\theta
Z)\}\Delta_{\chi+\ddot{\chi}}(Z)dZ
&=&\Delta_{\chi}(\theta^{-1})\frac{1}{\Gamma_{{\cal{P}}}(\chi)}
\int_{\cal{P}}\exp-\{\textrm{tr}Y\}\Delta_{\chi+\ddot{\chi}}(Y)dY\\
&=&\Delta_{\chi}(\theta^{-1}).
\end{eqnarray*}
\end{Pff}
\section{Riesz exponential families}In this section, we study the
natural exponential family generated by a Riesz measure. We first
review some basic concepts concerning exponential families and their
variance functions and introduce some notations.

For a positive measure on $\cal{A}$, we denote
$$\Theta(\mu)=\textrm{interior}\displaystyle\{\theta \in {\cal{A}}^{\ast}; \ L_{\mu}(\theta)<\infty\displaystyle\}$$
$$k_\mu=\log L_\mu$$
where $L_\mu$ and $k_\mu$ are respectively the Laplace transform and
the cumulant generating function of $\mu$.

The set $\mathcal{M}({\cal{A}})$ is now defined as the set of
positive measures $\mu$ such that $\mu$ is not concentrated on an
affine hyperplane of $\cal{A}$ and $\Theta(\mu)$ is not empty. For
$\mu$ in $\mathcal{M}({\cal{A}})$, the set of probability
$$F=F(\mu)=\{P(\theta,\ \mu)dX=\exp\{-\mbox{tr}(\theta X)-k_\mu\}\mu(dX); \ \theta \in \Theta(\mu)\}$$
\vskip0.1cm\n To each $\mu \in{\mathcal{M}}(E)$ and
$\theta\in\Theta(\mu)$, we associate the probability distribution on
${\cal{A}}$ \vskip0.2cm$\hfill
P(\theta,\mu)(dX)=exp\left(\langle\theta,X\rangle -k_\mu
(\theta)\right)\mu(dX).\hfill$\newline
 The set\vskip0.1cm$\hfill
 F=F(\mu)=\{P(\theta,\mu); \ \theta\in\Theta(\mu)\}\hfill$
 \vskip0.1cm\n is called the natural
exponential family (NEF) generated by $\mu$. We also say that $\mu$
is a basis of $F$. Note that a basis of $F$ is by no means unique.
If $\mu$ and $\nu$ are in ${\mathcal{M}}({\cal{A}})$, then it is
easy to check that $F(\mu)=F(\nu)$ if and only if there exist $a\in
{\cal{A}}$ and $b\in \reel$ such that $d\nu(X)=\exp(\langle
a,X\rangle +b)d\mu(X).$ Therefore, if $\mu$ is in
${\mathcal{M}}({\cal{A}})$ and $F=F(\mu)$, then
$${\mathcal{B}}_F=\{\nu\in{\mathcal{M}}({\cal{A}}); \ F(\nu)=F\}=
\{\exp(\langle a,X\rangle +b)\mu(dX); \ (a,b)\in {\cal{A}}\times
\reel\}$$ is the set of basis of $F.$
 \vskip0.2cm\n The function
$k_\mu$ is strictly convex and real analytic. Its first derivative
$k'_\mu$ defines a diffeomorphism between $\Theta(\mu)$ and its
image $M_F$. Since $k'_\mu (\theta)=\displaystyle\int X
P(\theta,\mu)(dX)$, $M_F$ is called the domain of the means of $F$.
The inverse function of $k'_\mu$ is denoted by $\psi_\mu$ and
setting $P(m,F)=P(\psi_\mu(m),\mu)$ the probability
 of $F$ with mean $m$, we have $$F=\left\{P(m,F);m\in M_F\right\},$$ which is the
parametrization of $F$ by the mean.\vskip0.1cm\n If $\mu$ and
$\nu=\exp(\langle a,X\rangle +b)\mu$ are two basis of $F$, then for
all $\theta\in D(\nu) =D(\mu)-a,$
\begin{equation}\label{N1}
k_\nu(\theta)=
 k_\mu(\theta+a)+b
 \end{equation}
 and for all $m\in M_F,$
 \begin{equation}\label{N2}
\psi_\nu(m)=\psi_\mu(m)-a. \end{equation} \n Now the covariance
operator of $P(m,F)$ is denoted by $V_F (m)$ and the map defined
from $M_F$ into $L_s({\cal{A}})$ by $m\longmapsto V_F
(m)=k''_\mu(\psi_\mu(m))$ is called the variance function of the NEF
$F$. It is easy proved that $V_F(m)=(\psi'_\mu(m))^{-1}$ and an
important feature of $V_F$ is that it characterizes $F$ in the
following sense: If $F$ and $F'$ are two NEFs such that $V_F (m)$
and $V_{F'}(m)$ coincide on a nonempty open set of $M_F \cap
M_{F'},$ then $F=F'$. In particular, knowledge of the variance
function gives knowledge of the NEF. \vskip0.2cm\n Let
$\varphi(X)=\delta(X)+\gamma$ be an affine transformation on
${\cal{A}}$, where $\delta\in GL({\cal{A}})$ and $\gamma\in
{\cal{A}}$ and let $F$ be some NEF, generated by $\mu$ we denote by
$\mu_{_{1}}=\varphi \ast \mu$ the image measure of $\mu$ by
$\varphi,$ then for all $\theta\in\Theta(\mu_{_{1}})=\delta^{\ast^{
-1}}(\Theta(\mu)),$
\begin{equation}\label{M1}
k_{\mu_{_{1}}}(\theta)=k_\mu(\delta^\ast(\theta))+\langle\theta,\gamma\rangle\cdot
 \end{equation} The following theorem gives a necessary and sufficient condition on ${\chi}$ so that
$R_{\chi}$ generates a natural exponential family.
\begin{theorem}\label{M(A)}
Let $\chi=\{\lambda_i, \ i \in I\}$ be in $\Xi$. Then the Riesz
measure $R_{\chi}$ is in $\mathcal{M}({\cal{A}})$ if and only if
$\lambda_{i}\not=0$, for $i \in \wp\cup{\cal{S}}$ .
\end{theorem}

\begin{Pff} \ $(\Leftarrow)$ Suppose that $\chi=\{\lambda_i, \ i \in I\}$
is in $\Xi$ such that $\lambda_{i}\not=0$ $\forall i \in
\wp\cup{\cal{S}}$. We have $\Theta(R_{\chi})$ is not empty since
it contains ${\cal{P}}^\ast$. We need to show that $R_{\chi}$ is
not concentrated on an affine hyperplane of $\cal{A}$. Write
$\chi=\displaystyle\sum_{i \in \wp\cup \mathcal{S}}\chi_{i}$ where
$\chi_{i} \in {\Xi}(i)$ (see \ref{Xi(i)}). Then
$R_{\chi}=\displaystyle\prod^{\ast}_{i \in \wp\cup
\mathcal{S}}R_{\chi_{i}}$, and it suffices to show that for any $i
\in \wp\cup{\cal{S}}$, $R_{\chi_{i}}$ is not concentrated on a
affine hyperplane of ${\cal{A}}_{i\preceq}$. In fact suppose that
there exists $i \in \wp\cup{\cal{S}}$ such that $R_{\chi_{i}}$ is
concentrated on a affine hyperplane $H$ of ${\cal{A}}_{i\preceq}$.
Then there exists $\psi \in \varepsilon^{i}$ such that
${\mathcal{T}}_{l}^{+}.e_{\psi}\subset H$. On the other hand,
there exist an element $a \in {\cal{A}}_{i\preceq}$ and an
hyperplane $H_0$ of ${\cal{A}}_{i\preceq}$ such that $H=a+H_0$.
Write
$${\cal{A}}_{i\preceq}=H_0+\reel X,$$
where $X \not \in H_0$. Let $(e_{lj}^k)_{1\leq k\leq n_{lj}}$ be a
basis of ${\mathcal{A}}_{lj}, \ (l,j) \in I_{i\preceq}\times
I_{i\preceq}$.

\n As ${\mathcal{A}}_{i\preceq}=\displaystyle\prod\limits _{(l,j)
\in I_{i\preceq}\times I_{i\preceq}}{\mathcal{A}}_{lj}$,  we can
write
$$X=\displaystyle\sum_{l , j\in I_{i\preceq}}\displaystyle\sum_{1\leq k\leq n_{lj}}\beta^k _{lj}e_{lj}^k,$$
where for $\ l , j\in I_{i\preceq}$ and $1\leq k\leq n_{lj}$,
$\beta^k _{lj}$, is a real number. Thus, we can write
$${\cal{A}}_{i\preceq}=H_0+\displaystyle\sum_{l , j\in
I}\displaystyle\sum_{1\leq k\leq n_{lj}}\reel e_{lj}^k.$$ Since the
dimension of $H_0$ is equal to $\dim({\cal{A}}_{i\preceq})-1$, then
there exist $l, j \in I_{i\preceq}$, and $1\leq k\leq n_{lj} $ such
that
$${\cal{A}}_{i\preceq}=H_0+\reel e_{lj}^k.$$
Let us consider the vectors
$$A_1=e_{lj}^k+\displaystyle\sum_{i\preceq j}\psi(j)e_{jj}^1 \hskip0.5cm\textrm{and} \hskip0.5cm\ A_2=2e_{lj}^k+\displaystyle\sum_{i\preceq j}\psi(j)e_{jj}^1\cdot$$
It is clear that $A_1$ and $A_2$ are in
${\mathcal{T}}_{l}^{+}.e_{\psi}\subset H$. Using the fact that for
$i \in \wp\cup{\cal{S}}$, $\lambda_i\not=0$, we have necessarily
$\psi(i)\not=0$, and we get $A_2-A_1=e_{lj}^k$ which is an element
of $H_0$. This is in contradiction with the fact that $X \not \in
H_0$. Thus for any $i \in \wp\cup{\cal{S}}$, $R_{\chi_{i}}$ is not
concentrated on a affine hyperplane of ${\cal{A}}_{i\preceq}$ and
Thereoem \ref{M(A)} is proved.

\n $(\Rightarrow)$ Suppose that $\lambda_i=0, \forall i \in \wp\cup
\mathcal{S}$. As the support of $R_{\chi}$ is
$\overline{{\cal{P}}}$, thus $R_{\chi}$ is not an element of
$\mathcal{M}({\cal{A}})$.
\end{Pff}

\n Next, we give the variance function of the Riesz exponential
family $F(R_{\chi})$ generated by $R_{\chi}$. For $X$ and $K$ in
${\cal{A}}$, we define the quadratic representation $P(X)$ by
$$P(X)K=X(KX)\cdot$$
It is symmetric, since we have $\langle P(X)K, \ L\rangle=\langle
P(X)L, \ K\rangle$.
\begin{theorem}\label{variance}
For any $m \in {\cal{P}}$,
\begin{eqnarray*}
V_{F(R_{\chi})}(m)=\displaystyle\sum_{i \in
I}\frac{1}{\lambda_i}(P(m_{i\preceq})-P(m_{i\prec}))=\displaystyle\sum_{i
\in I}\frac{1}{\lambda_i}P_i(P(m))\cdot
\end{eqnarray*}
\end{theorem}
For the proof of this theorem we were led to establish the following
intermediary result
\begin{lemma}\label{different X^-1}
The map
\begin{eqnarray*}
\varphi: X\rightarrow X^{-1}:  {\cal{P}}&\rightarrow&\reel\nonumber\\
Y&\mapsto&\textrm{tr}(X^{-1}Y)\nonumber
\end{eqnarray*}
is differentiable and its differential is
\begin{eqnarray*}
\varphi'(X)(K)=(X^{-1})'(K)=-P(X^{-1})(K),
\end{eqnarray*}
that is
\begin{eqnarray*}
(X^{-1})'(K)(Y)=-\textrm{tr}((X^{-1}(KX^{-1}))Y).
\end{eqnarray*}
\end{lemma}

\begin{Pff} \ We have
\begin{eqnarray*}
\displaystyle\lim_{t\rightarrow
0}\frac{1}{t}((X+tK)^{-1}-X^{-1})(Y)&=&\displaystyle\lim_{t\rightarrow
0}\frac{1}{t} \
\textrm{tr}((X+tK)^{-1}((X-(X+tK))X^{-1})Y)\\
&=&-\displaystyle\lim_{t\rightarrow 0}\frac{1}{t} \
\textrm{tr}((X+tK)^{-1}(KX^{-1})Y)\\
&=&-(X^{-1}(KX^{-1}))(Y)
\end{eqnarray*}
We now show that
\begin{eqnarray*}
(X+K)^{-1}-X^{-1}+X^{-1}(KX^{-1})=o(K)
\end{eqnarray*}
\small{\begin{eqnarray*}
\parallel(X+K)^{-1}-X^{-1}+X^{-1}(KX^{-1})\parallel&=&\parallel((X+K)^{-1}(X-(X+K)+(X+K)(X^{-1}K))X^{-1})\parallel\\
&=&\parallel((X+K)^{-1}(X-(X+K)+(X+K)(X^{-1}K))X^{-1})\parallel\\
&=&\parallel((X+K)^{-1}(K(X^{-1}K))X^{-1})\parallel\\
&\leq&
\parallel(X+K)^{-1}\parallel\parallel X^{-1}\parallel \parallel
K^2\parallel\cdot
\end{eqnarray*}}
Then  $$(X^{-1}(K))'=-X^{-1}(KX^{-1})=-P(X^{-1})(K).$$
\end{Pff}

\begin{Pff}{ \bf of Theorem \ref{variance}}
\ For $\theta \in {\cal{P}}^{\ast}$, we have
$L_{R_{\chi}}(\theta)=\Delta_{\chi}(\theta^{-1}).$ Then using
(\ref{genrateur power}), we get
\begin{eqnarray*}
k_{R_{\chi}}(\theta)=-\displaystyle\sum_{i \in
I}\lambda_{i}(\log\Delta_{\preceq i}(\theta)-\log\Delta_{\prec
i}(\theta))\cdot
\end{eqnarray*}
Therefore
\begin{eqnarray*}
k'_{R_{\chi}}(\theta)=-\displaystyle\sum_{i \in
I}\lambda_{i}((\theta_{\preceq i})^{-1}-(\theta_{\prec i})^{-1}),
\end{eqnarray*}
and from Lemma \ref{different X^-1}, we get
\begin{eqnarray*}
k^{''}_{R_{\chi}}(\theta)=\displaystyle\sum_{i \in
I}\lambda_{i}(P((\theta_{\preceq i})^{-1})-P((\theta_{\prec
i})^{-1})).
\end{eqnarray*}
It is easy to see that $(\theta_{\preceq
i})^{-1}=(\theta^{-1})_{\preceq i}$ and $(\theta_{\prec
i})^{-1}=(\theta^{-1})_{\prec i}$, then
\begin{eqnarray*}
k'_{R_{\chi}}(\theta)&=&-\displaystyle\sum_{i \in
I}\lambda_{i}((\theta^{-1})_{\preceq i}-(\theta^{-1})_{\prec
i})\cdot
\end{eqnarray*}
For $m \in {\cal{P}}$, such that $m=k'_{R_{\chi}}(\theta)$, we have
\begin{eqnarray*}
\theta^{-1}&=&-\displaystyle\sum_{i \in
I}\frac{1}{\lambda_{i}}(m_{\preceq i}-m_{\prec i}),
\end{eqnarray*}
then
\begin{eqnarray*}
(\theta^{-1})_{i\preceq}=-\frac{1}{\lambda_{i}}m_{i\preceq}, \ \
(\theta^{-1})_{i\prec}=-\frac{1}{\lambda_{i}}m_{i\prec}\cdot
\end{eqnarray*}
Therefore
\begin{eqnarray*}
V_{F(R_{\chi})}(m)=\displaystyle\sum_{i \in
I}\frac{1}{\lambda_i}(P(m_{i\preceq})-P(m_{i\prec}))\cdot
\end{eqnarray*}

\end{Pff}

\n Note that, in the particular case of the Wishart NEF, that is
when $\chi=\{\lambda, \forall i \in I\}$, we have

\begin{eqnarray*}
V_{F(R_{\chi})}(m)=\frac{1}{\lambda}\displaystyle\sum_{i \in
I}(P(m_{i\preceq})-P(m_{i\prec}))=\frac{1}{\lambda}P(m)\cdot
\end{eqnarray*}

\end{document}